\documentclass[11pt]{amsart}
\usepackage{amssymb, latexsym}
\usepackage{graphicx}
\theoremstyle{plain}
\newtheorem{theorem}{Theorem}
\newtheorem{corollary}{Corollary}

\theoremstyle{remark}

\newtheorem*{Remark 1}{Remark 1}
\newtheorem*{Remark 2}{Remark 2}
\newtheorem*{Remark 3}{Remark 3}
\newtheorem*{Remark 4}{Remark 4}

\numberwithin{equation}{section}

\begin{document}

\title[Cyclic to Random Shuffle Shuffles]%
 {Cyclic to Random Transposition Shuffles}

\author{Ross G. Pinsky}
\address{Department of Mathematics\\
Technion---Israel Institute of Technology\\
Haifa, 32000\\ Israel}
\email{ pinsky@math.technion.ac.il}
\urladdr{http://www.math.technion.ac.il/~pinsky/}
\thanks{}

\subjclass[2000]{60C05, 05A05, 05A15} \keywords{random shuffle, random permutation, total variation norm}
\date{}

\begin{abstract}
Consider  a permutation  $\sigma\in S_n$ as a deck of cards numbered from 1 to $n$ and laid out in a row,
 where
$\sigma_j$ denotes the number of the card that is in the $j$-th
position from the left.\rm\
We define  two cyclic to random transposition shuffles.
The first one works as follows: for $j=1,\cdots, n$, on the $j$-th step transpose the card that was \it originally\rm\
the $j$-th from the left  with a random card (possibly itself).
The second shuffle works as follows: on the $j$-th step, transpose the card that is \it currently\rm\ in the $j$-th
position from the left with
a random card (possibly itself).
For these shuffles, for  each $b\in[0,1]$, we calculate explicitly the limiting rescaled density function of $x,0\le x\le1$, for
the probability that a card with a number around $bn$ ends up in a position around $xn$,
and for each $x\in[0,1]$, we calculate the limiting rescaled density function of $b,0\le b\le 1$,  for
the probability that the card in a position around $xn$ will be a card with a number around $bn$.
These density functions all have a discontinuity at $x=b$, and for each  of them,
the supremum of the density is obtained by approaching the discontinuity from one side, and, for certain values
of the parameter,  the infimum of the density is
obtained by approaching the discontinuity from the other side.
\end{abstract}

\maketitle

\section{Introduction and Statement of Results}

Let $S_n$ denote  the symmetric group of permutations of $[n]\equiv\{1,\cdots, n\}$.
Our convention will be to  view  a permutation $\sigma\in S_n$ as a deck of cards numbered from 1 to $n$
and laid out in a row, where
$\sigma_j$ denotes the number of the card that is in the $j$-th position from the left.\rm\
In a recent paper \cite{P}, we analyzed the bias in the \it card cyclic to random insertion shuffle\rm:
remove and then randomly reinsert each of the $n$ cards exactly once,
the removal and reinsertion  being performed according to the \it original\rm\ left to right order
of the cards.
The novelty in this nonstandard shuffle  is that every card is removed and reinserted exactly once, unlike in any of the shuffles
one encounters in the literature.
The bias in this shuffle turned out to be surprisingly high, and possessed some interesting features.
We describe one of these features now, the one that is the impetus for the present article.

According to our convention, $\sigma^{-1}_j$ denotes the position occupied by  card number $j$.
Let $p_n(\text{id},\cdot)$ denote the probability distribution for the above shuffle when it starts from the identity
permutation.  Let $b\in[0,1]$ and let $\lim_{n\to\infty}b_n=b$, with $b_nn$ being an integer in $[n]$. It was shown that  the distribution
function $F_b(x)\equiv\lim_{n\to\infty}p_n(\text{id},\{\sigma^{-1}_{b_nn}\le xn\})$ for the limiting rescaled
position of a card with a number around $bn$ possesses a density $f_b(x), 0\le x\le 1$, which has a jump discontinuity; moreover,
$f_b$ is constant and at its minimum to the left of the discontinuity,  while as $x$ approaches
the point of discontinuity from the right, $f_b(x)$ approaches its supremum.
A similar phenomenon holds also for $h_x(b), 0\le b\le 1$, the density for the limiting rescaled card number in
a position around $xn$; more precisely, these densities also possess a discontinuity, the
supremum is approached as $b$ approaches the discontinuity from the left, while
 for a certain range of the parameter
value $x$, the infimum is  approached as $b$ approaches the point of discontinuity from the right.

In the present paper we consider the above quantities for two cyclic random transposition shuffles.
The first one works as follows: for $j=1,\cdots, n$, on the $j$-th step transpose the card that was \it originally\rm\
the $j$-th from the left  with a random card (possibly itself).
The method in this shuffle is simply  the method of the above card cyclic to random insertion shuffle transferred from
the context of insertions to  the context
of transpositions.
We call this shuffle the \it card cyclic to random transposition shuffle\rm.
The second shuffle works as follows: on the $j$-th step, transpose the card that is \it currently\rm\ in the $j$-th
position from the left with
a random card (possibly itself).
We call this shuffle  the \it position cyclic to random transposition shuffle\rm.

These two shuffles are a lot more similar to one another than the above card cyclic to random insertion shuffle
is to the corresponding \it position cyclic to random insertion shuffle\rm,
defined as follows: on step $j$, remove and randomly reinsert the card that is \it currently\rm\ in the $j$-th position from the left.
Indeed, it is easy to see that whereas by definition, every card gets removed and reinserted in the
card cyclic to random insertion shuffle, in general many cards do not get removed and reinserted at all
in the position cyclic to random insertion shuffle.
On the other hand, it is easy to see that in both the card cyclic to random transposition shuffle and the position
cyclic to random transposition shuffle, every card will  get transposed.

We will denote the probability measures on $S_n$ corresponding respectively to the
card cyclic to random transposition shuffle and the
position cyclic to random transposition shuffle starting from $\sigma\in S_n$ by $p_n^{\text{card}}(\sigma,\cdot)$
and  $p_n^{\text{pos}}(\sigma,\cdot)$.

\begin{theorem}\label{1}
\noindent i.
Under $p_n^{\text{card}}(\text{id},\cdot)$, the random variable $\sigma^{-1}_{b_nn}$, denoting the position of card
number $b_nn$, has the following behavior. Assume that $\lim_{n\to\infty}b_n=b\in[0,1]$ and that
$\lim_{n\to\infty}x_n=x\in[0,1]$. Then
$$
f^{\text{card}}_b(x)\equiv\lim_{n\to\infty}np_n^{\text{card}}(\text{id},\{\sigma^{-1}_{b_nn}=x_nn\})=
\begin{cases} e^{-1+b}+e^{-x}-e^{-1-x+b},\ x<b;\\
e^{-1+b}+e^{-x},\ x>b.
\end{cases}
$$
\noindent ii. Under $p_n^{\text{card}}(\text{id},\cdot)$, the random variable $\sigma_{x_nn}$, denoting the number of the card in position
$x_nn$, has the following behavior. Assume that $\lim_{n\to\infty}x_n=x\in[0,1]$ and that
$\lim_{n\to\infty}b_n=b\in[0,1]$.
Then
$$
h^{\text{card}}_x(b)\equiv\lim_{n\to\infty}np_n^{\text{card}}(\text{id},\{\sigma_{x_nn}=b_nn\})=
\begin{cases} e^{-1+b}+e^{-x},\ b<x;\\
e^{-1+b}+e^{-x}-e^{-1-x+b},\ b>x.\end{cases}
$$
\end{theorem}

\begin{theorem}\label{2}
\noindent i.
Under $p_n^{\text{pos}}(\text{id},\cdot)$, the random variable $\sigma^{-1}_{b_nn}$, denoting the position of card
number $b_nn$, has the following behavior. Assume that $\lim_{n\to\infty}b_n=b\in[0,1]$ and that
$\lim_{n\to\infty}x_n=x\in[0,1]$. Then
$$
f^{\text{pos}}_b(x)\equiv\lim_{n\to\infty}np_n^{\text{pos}}(\text{id},\{\sigma^{-1}_{b_nn}=x_nn\})=
\begin{cases}
 e^{-1+x}+e^{-b},\ x<b;\\
e^{-1+x}+e^{-b}-e^{-1-b+x},\ x>b.
\end{cases}
$$
\noindent ii. Under $p_n^{\text{pos}}(\text{id},\cdot)$, the random variable $\sigma_{x_nn}$, denoting the number of the card in position
$x_nn$, has the following behavior. Assume that $\lim_{n\to\infty}x_n=x\in[0,1]$ and that
$\lim_{n\to\infty}b_n=b\in[0,1]$.
Then
$$
h^{\text{pos}}_x(b)\equiv\lim_{n\to\infty}np_n^{\text{pos}}(\text{id},\{\sigma_{x_nn}=b_nn\})=
\begin{cases}
e^{-1+x}+e^{-b}-e^{-1-b+x},\ b<x;\\
e^{-1+x}+e^{-b},\ b>x.
\end{cases}
$$
\end{theorem}

\noindent \bf Remark 1.\rm\ Of course, part (ii) of each to the theorems follows imediately from part (i), and gives
 $h^\text{card}_x(b)=f^{\text{card}}_b(x)$ and
$h^\text{pos}_x(b)=f^{\text{pos}}_b(x)$.
Note also that it turns out that $f^{\text{pos}}_b(x)=h^{\text{card}}_x(b)$
and $h^{\text{pos}}_x(b)=f^{\text{card}}_b(x)$.

\begin{figure}
\includegraphics[scale=.7]{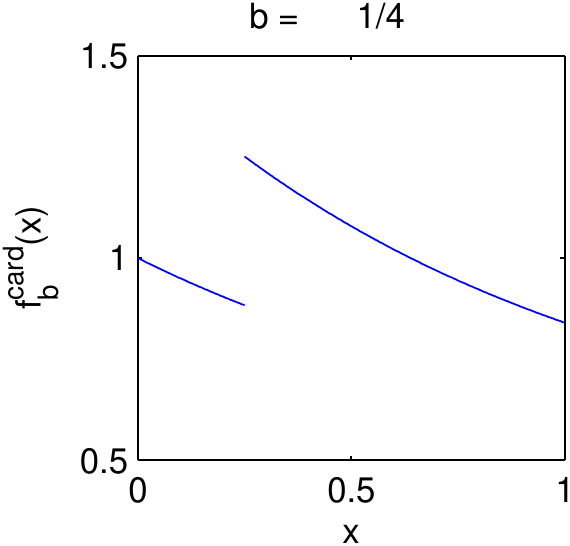}
\includegraphics[scale=.7]{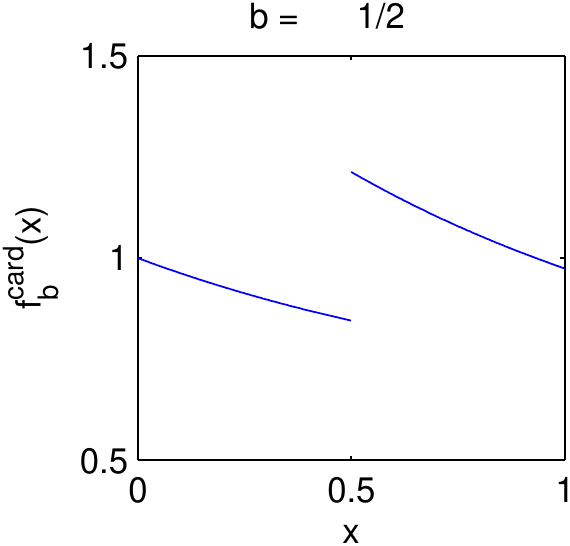}
\includegraphics[scale=.7]{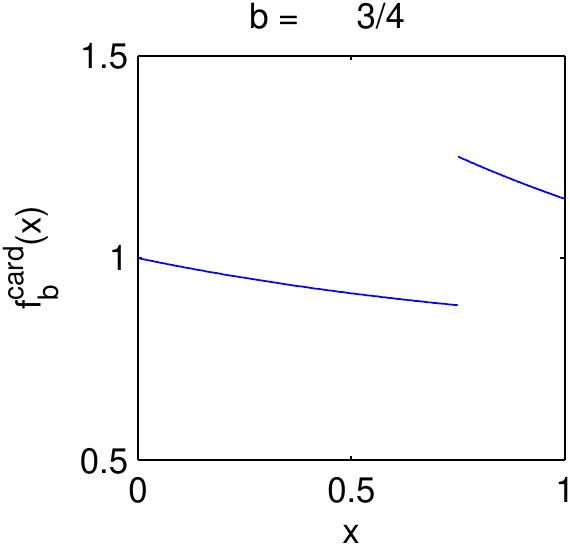}
\caption{Card cyclic to random transposition shuffle: density for limiting rescaled position of a card with a number around $bn$.}\label{F:f_b}
\end{figure}

\begin{figure}
\includegraphics[scale=.7]{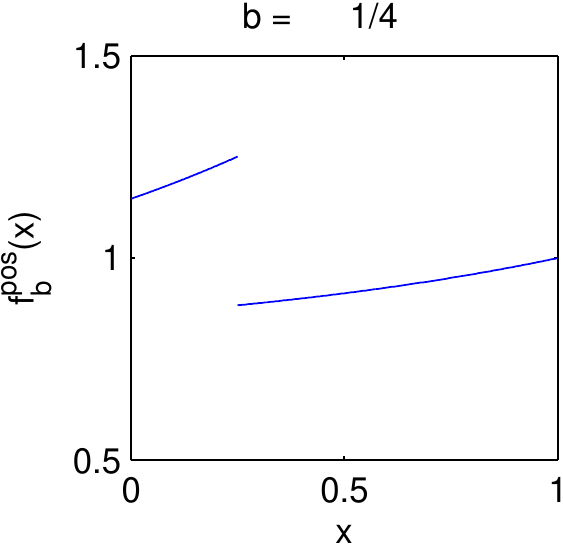}
\includegraphics[scale=.7]{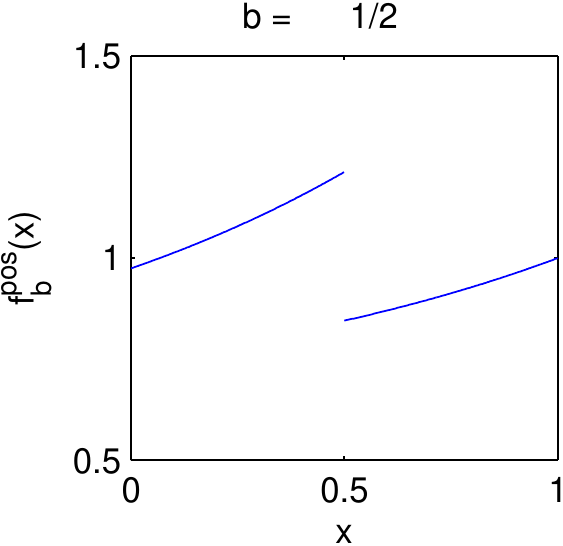}
\includegraphics[scale=.7]{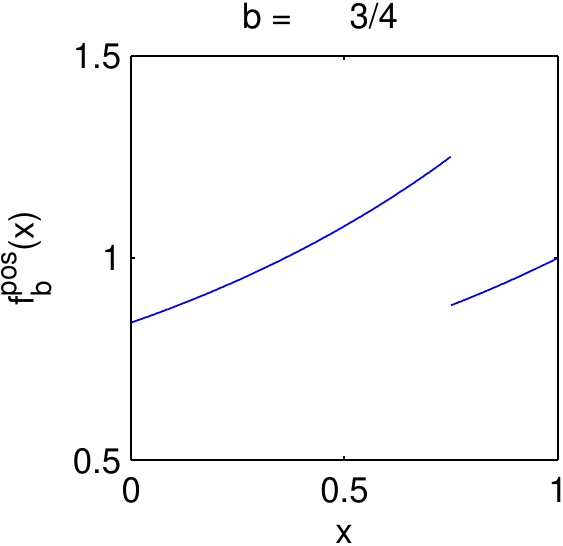}
\caption{Position cyclic to random transposition shuffle: density for limiting rescaled position of a card with a number  around $bn$.}\label{F:h_x}
\end{figure}

\bf\noindent Remark 2.\rm\  All four of the above densities have  a discontinuity when the variable  ($x$ or $b$, depending on the density)
is equal to the parameter ($b$ or $x$).
For the density $f^{\text{card}}_b(x)$,
 the supremum is approached  when the variable $x$   approaches the value of the parameter $b$
  from the right.  For the parameter $b$ in the range $[1-\log 2,1]\approx[.307,1]$,
 the infimum of the density  is approached when the variable approaches the value of the parameter from the left.
 For the parameter not in this range,
 the density  approaches a value less than 1 when the variable approaches the value of the parameter from the left, however
 the minimum of the density is attained at $x=1$.
For the density $f^{\text{pos}}_b(x)$,
the supremum is approached  when the variable $x$   approaches the value of the parameter $b$
from the left.  For the parameter $b$ in the range $[0,\log 2]\approx[0,.693]$,
the infimum of the density  is approached when the variable approaches the value of the parameter from the right.
For the parameter not in this range,
the density  approaches a value less than 1 when the variable approaches the value of the parameter from the right, however
the minimum of the density is attained at $x=0$.
Thus, under the card (position) cyclic to random transposition shuffle, the most likely position for a card to end up in
is a little to the right (left) of where it started, while the least likely position for it to end up in is either a little to the left
(right) of where it started or at the very end (beginning) of the deck, depending on where it started.
We conclude that the general
 phenomenon noted at the beginning of this paper with regard to the card cyclic to random insertion shuffle persists
 with the card cyclic and position cyclic to random transposition shuffles.
See figures \ref{F:f_b} and \ref{F:h_x} respectively for  graphs of $f^{\text{card}}_b$ and $f^{\text{pos}}_b$.
 The graphs of $h^{\text{card}}_x$  and $h^{\text{pos}}_x$ respectively are obtained from those of $f^{\text{pos}}_b$
 and $f^{\text{card}}_b$ by switching the labels of the $x$ and $b$ axes.

\medskip

Let $s$ and $t$ denote generically $b$ or $x$, and let $f_s^*(t)$ denote generically any of the above four density functions.
The supremum of $f_s^*$  is $e^{-1+s}+e^{-s}$, and the infimum is
$e^{-1+s}+e^{-s}-e^{-1}$. We have
$$
\sup_{0\le s\le 1}\sup_{0\le t\le 1}f_s^*(t)=1+e^{-1}\approx1.368,
$$
 with the supremum approached as
 $s$ and $t$ both approach 0 or 1, with one of the two variables strictly larger than the other one,
 the order depending on which of the four functions is considered. We have
 $$
 \inf_{0\le s\le 1}\inf_{0\le t\le 1}f_s^*(t)=2e^{-\frac12}-e^{-1}\approx.845,
$$
with the infimum approached
 as $s$ and $t$ both approach $\frac12$, with one of the two variables strictly larger than the other one,
 the order depending on which of the four functions is considered.

For the card cyclic to random transposition shuffle,
let
$$
E_{\text{pos}}^{\text{card}}(b)\equiv\int_0^1x f_b^{\text{card}}(x)dx\ \ \text{and}\ \
E_{\text{card}}^{\text{card}}(x)\equiv\int_0^1b h_x^{\text{card}}(b)db
$$
be respectively
the expected limiting rescaled position for a card with a number around $bn$ and
the expected limiting rescaled card number  to be found in a position around $xn$,
For
the position cyclic to random transposition shuffle,
let
$$
E_{\text{pos}}^{\text{pos}}(b)\equiv\int_0^1x f_b^{\text{pos}}(x)dx \ \ \text{and}\ \
E_{\text{card}}^{\text{pos}}(x)\equiv\int_0^1b h_x^{\text{pos}}(b)db
$$
be respectively
the expected limiting rescaled position for a card with a number around $bn$ and
the expected limiting
rescaled card number  to be found in a position around $xn$.

Direct calculations give the following corollary.
\begin{corollary}
\noindent i.
$$
E_{\text{pos}}^{\text{card}}(s)=E_{\text{card}}^{\text{pos}}(s)=1-\frac12e^{-1+s}-(1+s)e^{-1}.
$$
The maximum of this function occurs at $s=\log 2\approx .693$ with the value about .519,
and the minimum of this function occurs at $s=0$ with value about .448.

\noindent ii.
$$
E_{\text{card}}^{\text{card}}(s)=E_{\text{pos}}^{\text{pos}}(s)=\frac12e^{-s}+se^{-1}.
$$
The maximum of this function occurs at $s=1$ with the value about .552,
and the minimum of this function occurs at $s=1-\log 2\approx.307$ with value about .481.

\end{corollary}

 In the case of the card cyclic to random insertion shuffle, the density $f_b(x)$ degenerated when $b\to0$ to include
a $\delta$-mass at $x=0$ of weight $e^{-1}$. Using this as a starting point, we were able to show that
the total variation distance between the distribution $p_n(\text{id},\cdot)$ and the uniform distribution
converges to 1 as $n\to\infty$.
Recall that the total variation norm between two probability measures $\mu$ and $\nu$ on $S_n$ is defined
by
$$
||\mu-\nu||_{\text{TV}}=\sup_{A\subset S_n}(\mu(A)-\nu(A))=\frac12\sum_{\sigma\in S_n}|\mu(\sigma)-\nu(\sigma)|.
$$
  In the cases at hand, we don't know whether the total variation distance
between $p_n^{\text{card}}(\text{id},\cdot)$ and the uniform distribution or between
  $p_n^{\text{pos}}(\text{id},\cdot)$ and the uniform distribution goes to 1. Let $U_n$ denote the uniform distribution on $S_n$.
The above results just allow us to obtain the following extremely weak inequality.
Let $p^*_n(\text{id},\cdot)$ be generic notation
for either $p_n^{\text{pos}}(\text{id},\cdot)$ or $p_n^{\text{card}}(\text{id},\cdot)$.
Then
$$
\lim_{n\to\infty}||p^*_n(\text{id},\cdot)-U_n||_{\text{TV}}\ge\sup_{b\in[0,1]}\frac12\int_0^1|f_b^{\text{card}}(x)-1|dx\approx.08.
$$

As was the case with the card cyclic to random insertion shuffle, the card cyclic to random transposition shuffle does not
seem to have been studied before. On the other hand, the position cyclic to random transposition shuffle has been studied.
From the point of view of mixing times, it was studied in \cite{M} and \cite{MPS}, where it was shown that as $n\to\infty$,
 the number
of such shuffles needed to approach equilibrium is on the order $\log n$.
More in the spirit of this paper, the papers \cite{RB}, \cite{SS} and \cite{GM}
all studied various aspects of the distribution $p_n^{\text{pos}}(\text{id},\cdot)$, such as the limiting probability
of a derangement occurring, or the limiting expected number of fixed points.
In \cite{RB} it was shown that $p_n^{\text{pos}}(\text{id},\sigma)\ge\frac{2^{n-1}}{n^n}$, for all $\sigma\in S_n$,
and that the inequality is an equality
when $\sigma$ is  the permutation which cycles every card to the right.
In \cite{GM} it was also  shown that for $n\ge 18$, but not for $3\le n\le 17$, the identity permutation has the highest
probability; furthermore, as $n\to\infty$,
 $p_n^{\text{pos}}(\text{id},\text{id})\sim \frac{\frac1{\sqrt2}n^{\frac n2}e^{-\frac n2+\sqrt n-\frac14}}{n^n}$
 (\cite[p. 276]{RB}, \cite{GM}).
 By comparison, we note that in \cite{P} we showed that for the card cyclic to random insertion shuffle,
 one has the sharp inequalities  $\frac{2^{n-1}}{n^n}\le p_n(\text{id},\cdot)\le \frac{C_n}{n^n}$, where
 $C_n=\frac1{n+1}\binom {2n}n\sim\frac1{\sqrt\pi n^\frac32}\frac{4^n}{n^n}$ is the $n$-th Catalan number.
The results in \cite{RB}, \cite{SS}, and \cite{GM} do not at all allow one to determine whether or not
$\lim_{n\to\infty}||p^{\text{pos}}_n(\text{id},\cdot)-U_n||_{\text{TV}}=1$. It seems that the only result
in those papers that could be used to give an estimate on the total variation norm is the one
that says that under $p_n^{\text{pos}}(\text{id},\cdot)$, the  probability of a derangement
converges as $n\to\infty$ to  about .436 \cite{SS}, whereas under the uniform measure $U_n$ it converges to $e^{-1}\approx.367$, which gives
an upper bound on the total variation that is even weaker than the one we obtained above.

\noindent \bf Question:\rm\ Is it true that
$\lim_{n\to\infty}||p^{\text{card}}_n(\text{id},\cdot)-U_n||_{\text{TV}}=1$ and
$\lim_{n\to\infty}||p^{\text{pos}}_n(\text{id},\cdot)-U_n||_{\text{TV}}=1$?
\medskip

Theorem \ref{1} is proved in section 2 and Theorem \ref{2} is proved in section 3.

\section{Proof of Theorem \ref{1}}
Since $p_n^{\text{card}}(\text{id},\{\sigma^{-1}_j=a\})=p_n^{\text{card}}(\text{id},\{\sigma_a=j\})$,
part (ii) of the theorem follows immediately from part (i).
We now prove part (i).

For $a, j\in[n]$ with $a\neq j$, we consider $p_n^{\text{card}}(\text{id},\{\sigma^{-1}_j=a\})$, the probability
that card number $j$ ends up in position $a$.
The shuffle has $n$ steps, each of which is constituted by a transposition.
One way for $\{\sigma^{-1}_j=a\}$ to occur is for card number $j$ to move to position $a$ on the $j$-th step of the shuffle,
and then for it never to move again.
The probability of this occurring is $\frac1n(1-\frac1n)^{n-j}$.
The reader should convince himself that
if on the other hand,  card number $j$  moves to position $a$ on the $j$-th step of the shuffle, but
 is involved later on in another transposition, then it cannot end up in position $a$.
Another way for $\{\sigma^{-1}_j=a\}$  to occur is if $j<a$, card number $a$ is not moved before step $a$,
 on step $a$  card number $a$ is transposed with card number $j$,
and then card number $j$ is never transposed again. The probability of this is $\frac1n(1-\frac1n)^{n-1}$, if $j<a$.
If one the other hand, card number $j$ is transposed again after step $a$, then it can not end up in position $a$.

Assuming now that card number $j$ is not moved to position $a$ on the $j$-step or on the $a$-th step,  we consider how else one can end up with
the event $\{\sigma^{-1}_j=a\}$.
The reader should convince himself of the following facts. (Remember that all the following statements are being made under
the assumption that card number $j$ does not move to position $a$ on step $j$ or on step $a$.)
If card number $a$ gets transposed before the $a$-th step, then card number $j$ cannot
end up in  position $a$. (For example, say card number $a$ gets transposed for the first time on the $i$-th step, with $i<a$.
If card number $i$ was not transposed before the $i$-th step, then on the $i$-th step, card number $a$ was transposed with card number $i$.
So now after the $i$-th step, card number $i$ is  in position $a$. On every later step $k$, card number $k$ will
be transposed with a random card. Since we are assuming that card number $j$ was not moved to position $a$ on the $j$-th step,
there is no way for card number $j$ to end up in position $a$.
Similarly, if card number $i$ was transposed before step $i$, say at step $l$, and card number $l$ was not transposed before step $l$, then
on the $i$-th step, card number  $a$ was transposed with card number $l$. And the same logic as above shows that card number $j$ cannot
end up in position $a$.)

 If card number $a$ gets transposed for the first time at step $a$, and gets transposed with
a card whose number is less than or equal to $a$, but not equal to $j$,
then card number $j$ cannot end up in position $a$. (The reasoning is similar to the above reasoning.)

If card number $a$ gets transposed for the first time at step $a$, and gets transposed with a card whose number $l$ is
greater than $a$, but then card number $l$ gets transposed between step $a$ and step $l$, then for reasons similar to the above, card number $j$
cannot end up in position $a$.

However, if number card $a$ gets transposed for the first time at step $a$, and gets transposed with a card whose number $l$ is
greater than $a$, and then card number $l$ is not transposed between step $a$ and step $l$, then there
  is still a chance for card $j$ to end up in position $a$.
One way would be if
$l>j$ and if  on step $l$,  card number
$l$ is  transposed with card number $j$, and then after step $l$,  card number  $j$ is  never  transposed again (this last requirement
is possible because $l>j$).
Another way would be for card number
$l$ to get transposed on step $l$ with card number $m$ with $m>l$ and $m>j$, then for card number $m$ not
to get transposed between step $l$ and step $m$, then for card number $m$ to be transposed with  card number $j$ on step $m$,
and then for card number $j$ not to be transposed again after step $m$ (this last requirement is possible because
$m>j$). Continuing to argue in this vein, we arrive at the following formula.

Let $F_k$ denote the first step on which  card number $k$ is transposed. Of course, from the definition
of the shuffle one has $F_k\le k$.
Let  $T^i_k$ denote the last step strictly before  the $i$-th step on which card number $k$ was transposed.
If there is no such step then  define $T^i_k=\infty$. We will use the generic $P$ when considering probabilities
related to the random variables $F_k$ and $T^i_k$.
Then we have
\begin{equation}\label{cardinverse1}
\begin{aligned}
&p_n^{\text{card}}(\text{id},\{\sigma^{-1}_j=a\})=\frac1n(1-\frac1n)^{n-j}
+\frac1n(1-\frac1n)^{n-1}1_{j<a}+\\
&\sum_{m=1}^{n-a}\sum_{a\equiv i_0<i_1<\cdots <i_m\le n; i_m>j} P(F_a=a)\left(\prod_{l=1}^m P(T^{i_l}_{i_l}=i_{l-1}\right)P(T^{n+1}_j=i_m).
\end{aligned}
\end{equation}
We have $P(F_a=a)=(1-\frac1n)^{a-1}$, $P(T^{i_l}_{i_l}=i_{l-1})=\frac1n(1-\frac1n)^{i_l-i_{l-1}-1}$
and $P(T^{n+1}_j=i_m)=\frac1n(1-\frac1n)^{n-i_m}$, with this last equality holding because $i_m>j$.
Substituting this in \eqref{cardinverse1} gives
\begin{equation}\label{cardinverse2}
\begin{aligned}
&p_n^{\text{card}}(\text{id},\{\sigma^{-1}_j=a\})=\frac1n(1-\frac1n)^{n-j}+\frac1n(1-\frac1n)^{n-1}1_{j<a}+\\
&\sum_{m=1}^{n-a}\sum_{a\equiv i_0<i_1<\cdots <i_m\le n; i_m>j}
(1-\frac1n)^{n-m-1}(\frac1n)^{m+1}.
\end{aligned}
\end{equation}
If $m>j-a$, then the restriction $i_m>j$ on the inner sum above is superfluous and we have
$$
\sum_{a\equiv i_0<i_1<\cdots <i_m\le n; i_m>j}1=\binom{n-a}m.
$$
If $m\le j-a$, then we have
$$
\sum_{a\equiv i_0<i_1<\cdots <i_m\le n; i_m>j}1=\sum_{i_m=j+1}^n\binom{i_m-1-a}{m-1}.
$$
Using this with \eqref{cardinverse2} gives
\begin{equation}\label{cardinverse3}
\begin{aligned}
&p_n^{\text{card}}(\text{id},\{\sigma^{-1}_j=a\})=\frac1n(1-\frac1n)^{n-j}++\frac1n(1-\frac1n)^{n-1}1_{j<a}+\\
&\sum_{m=(j-a)^++1}^{n-a}(1-\frac1n)^{n-m-1}(\frac1n)^{m+1}\binom{n-a}m+\\
&\sum_{m=1}^{(j-a)^+}(1-\frac1n)^{n-m-1}(\frac1n)^{m+1}\left(\sum_{r=j+1}^n\binom{r-1-a}{m-1}\right),
\end{aligned}
\end{equation}
where any summation whose lower limit is greater than its upper limit is understood to vanish.

Let $X_{\text{Bin}(N,q)}$ denote a binomial random variable with parameters $N$ and $q$. We have
\begin{equation}\label{cardinverse4}
\sum_{m=(j-a)^++1}^{n-a}(1-\frac1n)^{n-m-1}(\frac1n)^{m+1}\binom{n-a}m=(1-\frac1n)^{a-1}\frac1nP(X_{\text{Bin}(n-a,\frac1n)}\ge (j-a)^++1).
\end{equation}
Similarly, for $r\in\{j+1,\cdots, n\}$,  we have
\begin{equation}\label{cardinverse5}
\begin{aligned}
&\sum_{m=1}^{(j-a)^+}(1-\frac1n)^{n-m-1}(\frac1n)^{m+1}\binom{r-1-a}{m-1}=\\
&\frac1{n^2}(1-\frac1n)^{n-r+a-1}P(X_{\text{Bin}(r-1-a,\frac1n)}\le (j-a)^+-1).
\end{aligned}
\end{equation}

Now let $a=x_n n$ and let $j=b_nn$, with $\lim_{n\to\infty}x_n=x$, $\lim_{n\to\infty}b_n=b$,
and $b\neq x$.
We conclude from \eqref{cardinverse4} that
\begin{equation}\label{cardinverse6}
\lim_{n\to\infty}n\sum_{m=(b_nn-x_nn)^++1}^{n-x_nn}(1-\frac1n)^{n-m-1}(\frac1n)^{m+1}\binom{n-x_nn}m=
 \begin{cases} e^{-x}(1-e^{-1+x}), \ \text{if}\ x>b;\\ 0,\ \text{if}\ x<b.\end{cases}
\end{equation}

For all $r\in\{b_nn+1,\cdots, n\}$, we have
\begin{equation}\label{cardinverse7}
\begin{aligned}
&\lim_{n\to\infty}P(X_{\text{Bin}(r-1-x_nn,\frac1n)}\le (b_nn-x_nn)^+-1)=1,\ \text{if}\
 x<b.
\end{aligned}
\end{equation}
It then follows from \eqref{cardinverse5} and \eqref{cardinverse7} that
\begin{equation}\label{cardinverse8}
\begin{aligned}
&\lim_{n\to\infty}n\sum_{m=1}^{(b_nn-x_nn)^+}(1-\frac1n)^{n-m-1}(\frac1n)^{m+1}\left(\sum_{r=b_nn+1}^n\binom{r-1-x_nn}{m-1}\right)=\\
&\lim_{n\to\infty}\sum_{r=b_nn+1}^n \frac1n(1-\frac1n)^{n-r+x_nn-1}=e^{-1-x}\int_b^1e^xdx=e^{-x}-e^{-1-x+b},\ \text{if}\ x<b.
\end{aligned}
\end{equation}
If $x>b$, then the left hand side of \eqref{cardinverse8} is identically 0 for large $n$.

It now  follows from \eqref{cardinverse3}, \eqref{cardinverse6}
and \eqref{cardinverse8} that
\begin{equation}\label{cardinverse9}
\lim_{n\to\infty}np_n^{\text{card}}(\text{id},\{\sigma^{-1}_{b_nn}=x_nn\})=
\begin{cases}
 e^{-1+b}+e^{-x},\ x>b;\\
e^{-1+b}+e^{-x}-e^{-1-x+b},\ x<b.
\end{cases}
\end{equation}
\hfill $\square$
\section{Proof of Theorem \ref{2}}
Since $p_n^{\text{pos}}(\text{id},\{\sigma^{-1}_j=a\})=p_n^{\text{pos}}(\text{id},\{\sigma_a=j\})$,
part (ii) of the theorem follows immediately from part (i).
We now prove part (i).

The analysis here is similar to that in the proof of Theorem \ref{1} so we will be less thorough here with  the explanations.
For $a, j\in[n]$ with $a\neq j$, we consider $p_n^{\text{pos}}(\text{id},\{\sigma^{-1}_j=a\})$, the probability
that card number $j$ ends up in position $a$.
The shuffle has $n$ steps, each of which is constituted by a transposition.
One way for $\{\sigma^{-1}_j=a\}$ to occur is for card number $j$ to be moved to position $a$ on the  $a$-th
step and then for card number $j$ never to move again. The probability of this is
$\frac1n(1-\frac1n)^{n-a}$.
Another way for $\{\sigma^{-1}_j=a\}$ to occur is if $a<j$,  card number $j$ does not  move before the $j$-th
step of the shuffle,   on the $j$-th step  it   moves to position $a$, and then  it never moves again
after the $j$-th step.
The probability of this is $\frac1n(1-\frac1n)^{n-1}$, if $a<j$. On the other hand,
if card number $j$ moves to position $a$ on the $a$-th step  and then moves again later on, it cannot end up in position $a$.
Similarly, if $a<j$, and  card number $j$ moves to position $a$ on the $j$-th step of the shuffle, and then moves
again later on, it cannot end up in position $a$.

Assuming now that card number $j$ is not moved to position $a$ on the $j$-step or on the $a$-th step,  we consider how else one can end up with
the event $\{\sigma^{-1}_j=a\}$.
By reasoning similar to that in the proof of Theorem \ref{1}, the only other way for this event to occur is
if there exists an $m\ge1$ and numbers $\{i_l\}_{l=1}^m$ satisfying
$j<i_1<\cdots<i_m\le n$, with $i_m>a$, such that card number $j$ does not move until step $j$, at step $j$ it is moved
to position $i_1$,
and then after that it moves again only at steps $i_l$, $l=1,\cdots, m$, with it moving at step $i_l$, $l=1, \cdots, m-1$,  to position $i_{l+1}$,
and with it moving at step $i_m$ to position $a$.
The probability of the above occurring for
a particular choice of $m$ and $\{i_l\}_{l=1}^m$ is
$$
\begin{aligned}
&(1-\frac1n)^{j-1}\frac1n(1-\frac1n)^{i_1-j-1}\times\cdots\times\frac1n(1-\frac1n)^{i_m-i_{m-1}-1}\frac1n(1-\frac1n)^{n-i_m}=\\
&(1-\frac1n)^{n-m-1}(\frac1n)^{m+1}.
\end{aligned}
$$
If $m>a-j$, then the restriction above that  $i_m>a$ is superfluous and we have
$$
\sum_{j<i_1<\cdots <i_m\le n; i_m>a}1=\binom{n-j}m.
$$
If $m\le a-j$, then we have
$$
\sum_{j<i_1<\cdots <i_m\le n; i_m>a}1=\sum_{i_m=a+1}^n\binom{i_m-1-j}{m-1}.
$$

From the above analysis we conclude that
\begin{equation}\label{posinverse1}
\begin{aligned}
& p_n^{\text{pos}}(\text{id},\{\sigma^{-1}_j=a\})=\frac1n(1-\frac1n)^{n-a}+\frac1n(1-\frac1n)^{n-1}1_{a<j}+\\
&\sum_{m=(a-j)^++1}^{n-j}(1-\frac1n)^{n-m-1}(\frac1n)^{m+1}\binom{n-j}m+\\
&\sum_{m=1}^{(a-j)^+}(1-\frac1n)^{n-m-1}(\frac1n)^{m+1}\left(\sum_{r=a+1}^n\binom{r-1-j}{m-1}\right),
\end{aligned}
\end{equation}
where any summation whose lower limit is greater than its upper limit is understood to vanish.
Noting that the right hand side of \eqref{posinverse1} is the right hand side of \eqref{cardinverse3} with the roles
of $a$ and $j$ switched, part (i) of the theorem follows from part (i) of Theorem \ref{1}.
\hfill $\square$.


\begin{thebibliography}{99}

\bibitem{GM}
 Goldstein, D. and  Moews, D.,  \emph{The identity is the most likely exchange shuffle for large $n$},  Aequationes Math. 65 (2003),  3-30.

\bibitem{M}
Mironov, I.,  \emph{(Not so) random shuffles of RC4}, Advances in cryptology—CRYPTO 2002, 304–319, Lecture Notes in Comput. Sci., 2442, Springer, Berlin, (2002),
304-319.

\bibitem{MPS} Mossel, E., Peres, Y. and Sinclair, A., \emph{Shuffling by semi-random transpositions},
 Foundations of Computer Science, 2004. Proceedings. 45th Annual IEEE Symposium, (2004), 572 - 581.

\bibitem{P} Pinsky, R.,  \emph{Probabilistic and  Combinatorial Aspects of the Card-Cyclic to Random Insertion Shuffle}, submitted.


\bibitem{RB} Robbins, D. P. and Bolker, E. D., \emph{The bias of three pseudorandom shuffles},
 Aequationes Math. 22 (1981), 268-292.

\bibitem{SS}
 Schmidt, F. and  Simion, R., \emph{ Card shuffling and a transformation on $S_n$},  Aequationes Math. 44 (1992), 11-34.

\end{thebibliography}
\end{document}